\documentclass{amsart}

\usepackage{amsmath, amsthm, amssymb, amsfonts, verbatim}

\def\h{\hat{h}}

\def\0{{\bf 0}}
\def\B{{\bf B}}

\def\f{{\mathfrak f}}
\def\p{{\mathfrak p}}

\def\Q{{\mathbb Q}}
\def\R{{\mathbb R}}
\def\Z{{\mathbb Z}}
\def\H{{\mathbb H}}
\def\C{{\mathbb C}}
\def\N{{\mathbb N}}
\def\H{{\mathbb H}}

\def\M{{\mathcal M}}

\def\br{\text{br}}

\def\tor{\text{tor}}
\def\ord{\text{ord}}

\theoremstyle{plain}

\newtheorem{thm}{Theorem}

\newtheorem{prop}[thm]{Proposition}
\newtheorem{lem}[thm]{Lemma}

\theoremstyle{definition}

\newtheorem*{theorem*}{Theorem}

\theoremstyle{remark}

\begin{document}

\title[Small Points on Elliptic Curves]{Small Rational Points on \\ Elliptic Curves Over Number Fields}

\date{August 3, 2005}

\author{Clayton Petsche}
\email{clayton@math.uga.edu}
\address{Department of Mathematics, 
         The University of Georgia, Athens, GA 30602-7403}

\begin{abstract}
Let $E/k$ be an elliptic curve over a number field.  We obtain some quantitative refinements of results of Hindry-Silverman, giving an upper bound for the number of $k$-rational torsion points, and a lower bound for the canonical height of non-torsion $k$-rational points, in terms of expressions depending explicitly on the degree $d=[k:\Q]$ of $k$ and the Szpiro ratio $\sigma$ of $E/k$.  The bounds exhibit only polynomial dependence on both $d$ and $\sigma$.
\end{abstract}

\maketitle


\begin{section}{Introduction}\label{Introduction}

Let $k$ be a number field of degree $d=[k:\Q]$, and let $E/k$ be an elliptic curve.  In \cite{Merel}, Merel used deep facts about the arithmetic of modular curves to prove that there is a universal bound $C(d)$ depending only on $d$ such that $|E(k)_{\tor}| \leq C(d)$.  Quantitative refinements of Merel's theorem due to Parent \cite{Parent} and Oesterl\'e (cf. \cite{HindrySilvermanII}) give bounds for $C(d)$ which depend exponentially on $d$, but it is still unknown whether one can bound $C(d)$ by an expression whose growth is polynomial in $d$.  Such a bound--or a proof that no such bound is possible--would be of value, both for its intrinsic interest and for its implications in cryptography; cf. \cite{ChengUchiyama}, \cite{KoblitzMenezes}.

In this paper we give an explicit polynomial bound on $|E(k)_{\tor}|$ depending on $d$ and the Szpiro ratio $\sigma$, a certain quantity associated to the elliptic curve $E/k$, which we will now define.  Recall that the conductor $\f_{E/k}$ and the minimal discriminant $\Delta_{E/k}$ are certain integral ideals of ${\mathcal O}_k$ that are supported on the primes at which $E/k$ has bad reduction.  The Szpiro ratio is given by
\begin{equation}\label{szpiroratio}
\sigma = \frac{\log |\N_{k/\Q}(\Delta_{E/k})|}{\log |\N_{k/\Q}(\f_{E/k})|}
\end{equation}
when $E/k$ has at least one place of bad reduction; by convention we put $\sigma=1$ if $E/k$ has everywhere good reduction.  

\begin{thm}\label{krationaltor}
Let $k$ be a number field of degree $d=[k:\Q]$, and let $E/k$ be an elliptic curve with Szpiro ratio $\sigma$.  Then
\begin{equation}\label{torsionbound}
|E(k)_{\tor}| \leq  c_1d\sigma^2\log(c_2d\sigma^2),
\end{equation}
where $c_1=134861$ and $c_2=104613$.  
\end{thm}

Strictly speaking $(\ref{torsionbound})$ is not a {\em uniform} bound, in that the right-hand-side depends on the elliptic curve $E/k$.  However, the result is perhaps more interesting in view of the fact that one generally expects $\sigma$ to be small.  More precisely, let $\Sigma(k)$ denote the set (with multiplicities) of Szpiro ratios $\sigma$ of all elliptic curves $E/k$.  A well-known conjecture of Szpiro \cite{Szpiro} asserts that $\Sigma(k)$ is bounded, and that $6$ is its largest limit point.  (As pointed out by Masser \cite{Masser}, $6$ is in fact a limit point.)  An analogue of Szpiro's conjecture is known to hold in the function field case, and in the number field case it can be shown to follow from the ABC conjecture.  

Note that if Szpiro's conjecture is true, then Theorem~\ref{krationaltor} gives a uniform bound on $|E(k)_{\tor}|$ in terms of $k$ only.  The first result showing that Szpiro's conjecture implies such a uniform bound is due to Frey; the argument was written up by Flexor-Oesterl\'e \cite{FlexorOesterle}.  Also in \cite{FlexorOesterle}, the authors treated the special case of everywhere good reduction ($\sigma=1$), giving a bound on $|E(k)_{\tor}|$ that is exponential in $d$.  In \cite{HindrySilvermanII} Hindry-Silverman improved this to a bound of $O(d\log d)$; thus our bound $(\ref{torsionbound})$ recovers theirs (with slightly different constants) in this special case.  In the general case, Hindry-Silverman \cite{HindrySilvermanI} have given a bound on $|E(k)_{\tor}|$ that is exponential in both $d$ and $\sigma$.  Thus the main interest in our bound $(\ref{torsionbound})$ is in its explicit nature, and in the fact that it exhibits only polynomial growth in both $d$ and $\sigma$.

Finally, let us denote by $\Sigma_d$ the union of the sets $\Sigma(k)$ over all number fields $k$ of degree $d=[k:\Q]$.  In view of Merel's Theorem and Theorem~\ref{krationaltor}, it is natural to ask whether $\sup\Sigma_d$ and $\limsup\Sigma_d$ are finite, and if so, how they depend on $d$.  Upper bounds on these quantities would amount to a significant strengthening of Szpiro's conjecture, however, and would therefore seem to lie well beyond the scope of present techniques.

\medskip

A problem somewhat related to counting rational torsion points is that of giving a lower bound on the N\'eron-Tate canonical height $\h(P)$, for non-torsion rational points $P$, depending explicitly on the relevant data associated to the elliptic curve $E/k$.  In particular, a conjecture of Lang asserts that if $P\in E(k)$ is not a torsion point, then 
\begin{equation}\label{langconj}
\h(P)\geq c\log|\N_{k/\Q}(\Delta_{E/k})|,
\end{equation}
where $c=c(k)>0$ is a constant depending only on $k$.  We show that one can take for $c$ a certain expression depending explicitly on $d$ and $\sigma$.

\begin{thm}\label{langs}
Let $k$ be a number field of degree $d=[k:\Q]$, and let $E/k$ be an elliptic curve with minimal discriminant $\Delta_{E/k}$ and Szpiro ratio $\sigma$.  Then 
\begin{equation}\label{kratlangs}
\h(P)\geq c(d,\sigma)\log|\N_{k/\Q}(\Delta_{E/k})|
\end{equation}
for all non-torsion points $P\in E(k)$, where 
\begin{equation}\label{cdsigma}
c(d,\sigma)=\frac{1}{10^{15}d^3\sigma^6\log^2(c_2d\sigma^2)},
\end{equation}
and $c_2=104613$.  
\end{thm}

A consequence of Theorem~\ref{langs} is that Szpiro's conjecture implies Lang's conjecture; this fact was originally proved by Hindry-Silverman \cite{HindrySilvermanI}, who showed that $(\ref{langconj})$ holds with a value of $c$ depending exponentially on $d$ and $\sigma$.   Thus again the main interest in $(\ref{kratlangs})$ is in the fact that $c(d,\sigma)$ exhibits only polynomial decay in $d$ and $\sigma$.  Compare with the results of David \cite{David}, who uses methods from transcendence theory to obtain a similar bound.  

\medskip

Let us now briefly summarize our approach, in which we extend the methods of Hindry-Silverman \cite{HindrySilvermanII} to include a treatment of the places of bad reduction.  Denote by $\M_k$ the set of all places of $k$, and given $v\in\M_k$, let $k_v$ be the completion of $k$ at $v$.  Let $\h: E(\bar{k})\to[0,+\infty)$ denote the N\'eron-Tate canonical height, and recall that given a point $P\in E(k)\setminus\{O\}$, we have the local decomposition $\h(P) = \sum_{v\in\M_k}\frac{d_v}{d}\lambda_v(P)$, where $d_v=[k_v:\Q_v]$ is the local degree and $\lambda_v:E(k_v)\setminus\{O\}\to\R$ is the appropriately normalized N\'eron local height function (cf. \cite{SilvermanII}, $\S$VI.1.)  Given a set $Z=\{P_1,\dots,P_N\}\subset E(k)$ of $N$ distinct small rational points, we estimate the height-discriminant sum 
\begin{equation}\label{heightdisc}
\Lambda(Z) = \frac{1}{N^2}\sum_{1\leq i,j\leq N}\h(P_i-P_j),
\end{equation}
from above globally via the parallelogram law, and from below locally using the decomposition $\Lambda(Z)=\sum_{v\in\M_k}\frac{d_v}{d}\Lambda_v(Z)$, where
\begin{equation}\label{heightdiscloc}
\Lambda_v(Z) = \frac{1}{N^2}\sum_{\stackrel{1\leq i,j\leq N}{i\neq j}}\lambda_v(P_i-P_j).
\end{equation}

In order to obtain the necessary archimedean estimates we follow \cite{HindrySilvermanII}, using the pigeonhole principle to pass to a subset of $Z$ of positive density, all of whose points are close to each other in $E(k_v)$ at a particular archimedean place $v$.  At the non-archimedean places we give an analytic lower bound on $\Lambda_v(Z)$ in terms of the valuation of the minimal discriminant (cf. Lemma~\ref{lambdalemma} below).  Although we will not specifically require this interpretation in the present paper, this inequality can be viewed as a quantitative form of a local equidistribution principle for small points, as developed in \cite{BakerPetsche}.  Finally, it is worth noting that we do not decompose the minimal discriminant into ``small power'' and ``large power'' parts as Hindry-Silverman do in \cite{HindrySilvermanI}.  Instead, we assemble the local information at the different places of bad reduction into global information by a simple application of Jensen's inequality.

The author would like to express his gratitude to Matt Baker for his many valuable suggestions concerning this paper.

\end{section}


\begin{section}{Elliptic Curves over Number Fields}
In this section we fix some notation, and we review some of the relevant facts concerning elliptic curves over number fields.

Let ${\mathcal O}_k$ be the ring of integers of a number field $k$, and denote by $\M_k^\infty$ and $\M_k^0$ the set of all archimedean and non-archimedean places of $k$ respectively.  For each place $v\in\M_k$ we will select the normalized absolute value $|\text{ }|_v$ that, when restricted to $\Q$, coincides with one of the usual archimedean or $p$-adic absolute values.  Given a place $v\in\M_k^0$ lying over the rational prime $p$, let ${\mathcal O}_v$ and $\M_v$ be the ring of integers and maximal ideal in $k_v$ respectively.  Let $\p_v={\mathcal O}_k\cap\M_v$ denote the prime ideal of ${\mathcal O}_k$ corresponding to the place $v$, and let $\pi_v\in\M_v$ be a uniformizer.  Thus $|\N_{k/\Q}(\p_v)|=|{\mathcal O}_v/\M_v|=p^{f_v}$ and $|\pi_v|_v=p^{-1/e_v}$, where $f_v$ and $e_v$ are the residual degree and ramification index, respectively.  Recall that $d_v=[k_v:\Q_p]=e_vf_v$.  

Let $E/k$ be an elliptic curve, and let $j_E$ denote its $j$-invariant.  The conductor and minimal discriminant of $E/k$ are the integral ideals 
\begin{equation}\label{conddiscdef}
\f_{E/k}  = \prod_{v\in\M_k^0}\p_v^{\eta_v}, \text{ and }  \Delta_{E/k}  = \prod_{v\in\M_k^0}\p_v^{\delta_v}
\end{equation}
of ${\mathcal O}_k$, respectively, where the exponents $\eta_v$ and $\delta_v$ are given as follows.  Fix a place $v\in\M_k^{0}$, and define $\delta_v=\ord_v(\Delta_v)$, where $\Delta_v\in{\mathcal O}_v$ denotes the discriminant of a minimal Weierstrass equation for $E/k_v$; thus $|\Delta_v|_v=|\pi_v|_v^{\delta_v}$.  The exponent $\eta_v$ of the conductor has a rather complicated Galois-theoretic definition (cf. \cite{SilvermanII}, $\S$ IV.10); alternatively, it can be characterized by Ogg's formula
\begin{equation}\label{ogg}
\delta_v = \eta_v +m_v -1,
\end{equation}
where $m_v$ is the number of components on the special fiber of the minimal proper regular model of $E/k_v$ (cf. \cite{SilvermanII}, $\S$ IV.11).  An immediate consequence of $(\ref{ogg})$ is that $\eta_v\leq\delta_v$, from which we deduce the lower bound $\sigma\geq1$ on the Szpiro ratio. 

Let $E_0(k_v)$ denote the subgroup of $E(k_v)$ consisting of those points whose reduction (with respect to a minimal Weierstrass equation for $E/k_v$) is non-singular, and let
\begin{equation}\label{cvdef}
c_v  = |E(k_v)/E_0(k_v)|
\end{equation}
be the cardinality of the quotient.  It is well known that the set $E_0(k_v)$, and thus also the number $c_v$, do not depend on the choice of the minimal Weierstrass equation used to define them.

Finally, we recall that the various data discussed above are governed by the reduction type of $E/k_v$.  To be precise, if $E/k_v$ has good reduction then $\delta_v=\eta_v=0$ and $c_v=1$; if $E/k_v$ has split multiplicative reduction then $\eta_v=1$ and $c_v=\delta_v$; and finally, in all other cases we have $c_v\leq4$.  For proofs of these assertions cf. \cite{SilvermanII}, $\S$ IV.9.

\medskip

Following Rumely \cite{Rumely}, it will be useful to decompose the non-archimedean N\'eron local height function $\lambda_v:E(k_v)\setminus\{O\}\to\R$ into a sum
\begin{equation}\label{decomp}
\lambda_v(P-Q) = i_v(P,Q)+j_v(P,Q),
\end{equation}
where $i_v$ is a nonnegative arithmetic intersection term, and $j_v$ is bounded (cf. also \cite{BakerPetsche}).  These functions are most naturally described by considering the cases of integral and non-integral $j$-invariant separately:

{\em Case 1:} $|j_E|_v\leq1$.  Then $j_v$ is identically zero and $\lambda_v(P-Q) = i_v(P,Q)$.

{\em Case 2:} $|j_E|_v>1$.  By Tate's uniformization theory \cite{Tate} we have maps
\begin{equation}\label{tate}
E(k_v) \stackrel{\sim}{\longrightarrow} k_v^\times/q^\Z \to \R/\Z,
\end{equation}
where $q\in k_v^\times$ with $|q|_v=|1/j_E|_v<1$; here the first map is the Tate isomorphism, and the second map is given by $u\mapsto\log|u|_v/\log|q|_v$.  Let $r$ denote the composition of the two maps in $(\ref{tate})$.  Then 
\begin{equation*}
j_v(P,Q) = \frac{1}{2}\B_2(r(P-Q))\log|j_E|_v,
\end{equation*}
where 
\begin{equation*}
\B_2(t) = (t-[t])^2 - \frac{1}{2}(t-[t]) + \frac{1}{12} = \frac{1}{2\pi^2}\sum_{m\in\Z\setminus\{0\}}\frac{1}{m^2}e^{2\pi imt}
\end{equation*}
is the periodic second Bernoulli polynomial.  If $r(P)\neq r(Q)$ then $i_v(P,Q)=0$.  On the other hand if $r(P)=r(Q)$, then we can select representatives $u(P)$, $u(Q)\in k_v^\times/q^\Z$ for $P$ and $Q$ under the Tate isomorphism such that $|u(P)|_v=|u(Q)|_v$; in this case
\begin{equation*}
i_v(P,Q) = -\log|1-u(P)/u(Q)|_v.
\end{equation*}
Finally, we note for future reference that $E_0(k_v)=\ker(r)$; in other words, a point $P\in E(k_v)$ has non-singular reduction if and only if $r(P)=0$.

\end{section}


\begin{section}{A Few Preliminary Lemmas}
The following lower bound on the non-archimedean local sum $(\ref{heightdiscloc})$ is a variant of \cite{HindrySilvermanIII}, Prop. 1.2.  To be precise, our bound $(\ref{lambdalowerbound})$ is given in terms of the valuation of the minimal discriminant, rather than the valuation of the $j$-invariant as in \cite{HindrySilvermanIII}, a distinction that is only relevant at the places of additive reduction.

\begin{lem}\label{lambdalemma}
Let $E/k$ be an elliptic curve, let $v$ be a non-archimedean place of $k$, and let $\Delta_v\in {\mathcal O}_v$ be the discriminant of a minimal Weierstrass equation for $E/k_v$.  If $Z\subset E(k_v)$ is a set of N distinct $k_v$-rational points, then
\begin{equation}\label{lambdalowerbound}
\Lambda_v(Z) \geq \Big(\frac{1}{c_v^2}-\frac{1}{N}\Big)\frac{1}{12}\log|1/\Delta_v|_v.
\end{equation}
\end{lem}
 
\begin{proof}
First, note that
\begin{equation}\label{jdiscineq}
\log^+|j_E|_v\leq\log|1/\Delta_v|_v,
\end{equation}
with equality if and only if $E/k_v$ has good or multiplicative reduction (cf. \cite{SilvermanI}, Prop. VII.5.1).  Also, the N\'eron local height function satisfies the lower bound $\lambda_v(P)\geq\frac{1}{12}\log|1/\Delta_v|_v$ for points $P\in E_0(k_v)\setminus\{O\}$ (cf. \cite{SilvermanII}, Thm. VI.4.1).  It follows from this and the decomposition $(\ref{decomp})$ that
\begin{equation}
i_v(P,Q) \geq \frac{1}{12}(\log|1/\Delta_v|_v - \log^+|j_E|_v)\geq0
\end{equation}
for $P-Q\in E_0(k_v)\setminus\{O\}$.  As $i_v$ is nonnegative, it follows that
\begin{equation}\label{m1}
\sum_{\stackrel{1\leq i,j\leq N}{i\neq j}}i_v(P_i,P_j) \geq (M-N)\frac{1}{12}(\log|1/\Delta_v|_v - \log^+|j_E|_v),  
\end{equation}
where $M$ is the number of ordered pairs $(P_i,P_j)\in Z\times Z$ with $P_i-P_j\in E_0(k_v)$.  For each coset $C\in E(k_v)/E_0(k_v)$, let $N_C=|Z\cap C|$.  Then
\begin{equation}\label{m2}
\begin{split}
M &  = \sum_{\stackrel{1\leq i,j\leq N}{P_i-P_j\in E_0(k_v)}}1 \\
	& = \sum_{C\in E(k_v)/E_0(k_v)}N_C^2 \\
	& = \sum_{C\in E(k_v)/E_0(k_v)}\Big\{\Big(N_C-\frac{N}{c_v}\Big)^2+\frac{2N_C N}{c_v}-\frac{N^2}{c_v^2}\Big\} \\
	& \geq \frac{2N}{c_v}\sum_{C\in E(k_v)/E_0(k_v)}N_C-\frac{N^2}{c_v^2} \sum_{C\in E(k_v)/E_0(k_v)}1 \\
	& = \frac{N^2}{c_v}.
\end{split}
\end{equation}
Combining $(\ref{m1})$ and $(\ref{m2})$ we find that
\begin{equation}\label{m3}
\frac{1}{N^2}\sum_{\stackrel{1\leq i,j\leq N}{i\neq j}}i_v(P_i,P_j) \geq \Big(\frac{1}{c_v}-\frac{1}{N}\Big)\frac{1}{12}(\log|1/\Delta_v|_v - \log^+|j_E|_v).
\end{equation}

We will now show that
\begin{equation}\label{jlowerbound}
\frac{1}{N^2}\sum_{\stackrel{1\leq i,j\leq N}{i\neq j}}j_v(P_i,P_j) \geq \Big(\frac{1}{c_v^2}-\frac{1}{N}\Big)\frac{1}{12}\log^+|j_E|_v.
\end{equation}
For if $|j_E|_v\leq1$ then both sides of $(\ref{jlowerbound})$ are zero.  On the other hand, assume that $|j_E|_v>1$.  Since the map $r$ is trivial on the identity component $E_0(k_v)$, it is well-defined on the quotient $E(k_v)/E_0(k_v)$, which has cardinality $c_v$.  Therefore $r(E(k_v))\subseteq\langle1/c_v\rangle\subset\R/\Z$, and we conclude that $e^{2\pi imr(P_j)}=1$ for all $P_j\in Z\subset E(k_v)$, whenever $c_v\mid m$.  It follows that
\begin{equation}
\begin{split}
\sum_{\stackrel{1\leq i,j\leq N}{i\neq j}}j_v(P_i,P_j) & = \sum_{1\leq i,j\leq N}j_v(P_i,P_j) -\frac{N}{12}\log|j_E|_v \\
	& = \frac{\log|j_E|_v}{4\pi^2}\sum_{m\in\Z\setminus\{0\}}\frac{1}{m^2}\Big|\sum_{1\leq j\leq N}e^{2\pi imr(P_j)}\Big|^2 -\frac{N}{12}\log|j_E|_v \\
	& \geq \frac{N^2\log|j_E|_v}{4\pi^2}\sum_{m\in c_v\Z\setminus\{0\}}\frac{1}{m^2} -\frac{N}{12}\log|j_E|_v \\
	& = \frac{N^2\log|j_E|_v}{12c_v^2} -\frac{N}{12}\log|j_E|_v,
\end{split}
\end{equation}
which implies $(\ref{jlowerbound})$.

Finally, by the decomposition $(\ref{decomp})$ and the lower bounds $(\ref{m3})$ and $(\ref{jlowerbound})$, and using the inequality $(\ref{jdiscineq})$, we have
\begin{equation*}
\begin{split}
\Lambda_v(Z) & = \frac{1}{N^2}\sum_{\stackrel{1\leq i,j\leq N}{i\neq j}}(i_v(P_i,P_j)+j_v(P_i,P_j)) \\
	& \geq \Big(\frac{1}{c_v^2}-\frac{1}{c_v}\Big)\frac{1}{12}\log^+|j_E|_v + \Big(\frac{1}{c_v}-\frac{1}{N}\Big)\frac{1}{12}\log|1/\Delta_v|_v \\
	& \geq \Big(\frac{1}{c_v^2}-\frac{1}{c_v}\Big)\frac{1}{12}\log|1/\Delta_v|_v + \Big(\frac{1}{c_v}-\frac{1}{N}\Big)\frac{1}{12}\log|1/\Delta_v|_v \\
	& = \Big(\frac{1}{c_v^2}-\frac{1}{N}\Big)\frac{1}{12}\log|1/\Delta_v|_v,
\end{split}
\end{equation*}
which completes the proof.
\end{proof}

The following archimedean analogue of Lemma~\ref{lambdalemma} is a quantitative refinement of a result due to Elkies (cf. \cite{LangAT}, $\S$VI.); for a detailed proof see \cite{BakerPetsche}, Appendix A. 

\begin{lem}\label{lambdalemmaarch}
Let $E/k$ be an elliptic curve with j-invariant $j_E$, and let $Z\subset E(k)$ be a set of N distinct $k$-rational points.  If $v$ is an archimedean place of $k$, then
\begin{equation}\label{elkiesbound}
\Lambda_v(Z) \geq -\frac{\log N}{2N} - \frac{1}{12N}\log^+|j_E|_{v} - \frac{16}{5N}.
\end{equation}
\end{lem} 

In order to obtain the necessary non-archimedean estimates we will require the following lemma due to Hindry-Silverman (cf. Lemma 1 of \cite{HindrySilvermanII} and Prop. 2.3 of \cite{HindrySilvermanI}.)  Let $j:\H=\{\tau\in\C\,\mid\,\Im(\tau)>0\}\to\C$ denote the modular j-function, and let $L=\Z+\tau\Z$ be a normalized lattice with $\tau\in\H$; thus $j(\tau)$ is the j-invariant of the elliptic curve $\C/L$.  Let $\lambda:(\C/L)\setminus\{0\}\to\R$ denote the N\'eron function, as given in \cite{SilvermanII}, $\S$ VI.3.

\begin{lem}[Hindry-Silverman]\label{hind}
If $z=r_1+r_2\tau\in\C\setminus\{0\}$, where $r_1,r_2\in\R$ and $\max\{|r_1|, |r_2|\}\leq1/24$, then
\begin{equation*}
\lambda(z) \geq \frac{1}{288}\max\{1,\log|j(\tau)|\}.
\end{equation*}
\end{lem}

Finally, we will require the following inequality, the proof of which is elementary and given in \cite{BakerPetsche}, $\S$ 6.
\begin{lem}\label{tlem}
Let $N\geq1$ satisfy the bound $N\leq A\log N+B$ for constants $A>0,B\geq0$.  
Then $N\leq (\frac{e}{e-1})(A\log A+B)$. 
\end{lem}
\end{section}


\begin{section}{A Bound on the Number of Small Points}

In this section we will prove the results stated in the introduction, which are both consequences of the following bound on the number of small rational points.

\begin{prop}\label{smallprop}
Let $k$ a number field of degree $d=[k:\Q]$, and let $E/k$ be an elliptic curve with Szpiro ratio $\sigma$.  Then
\begin{equation}\label{kratmain}
\Big|\Big\{P\in E(k)\,\Big|\, \h(P) \leq \frac{\log |\N_{k/\Q}(\Delta_{E/k})|}{2^{13}3d\sigma^2}\Big\}\Big| \leq c_1d\sigma^2\log(c_2d\sigma^2),
\end{equation} 
where $c_1=134861$ and $c_2=104613$.  
\end{prop}

\begin{proof}
To ease notation we will henceforth suppress the subscripts on the notations $\N_{k/\Q}$, $\Delta_{E/k}$, and $\f_{E/k}$.  Let $S$ denote the set on the left-hand-side of $(\ref{kratmain})$, and let $N$ denote the largest integer satisfying $|S|\geq 24^2(N-1)+1$; thus $|S|\leq 24^2N$.  We will show that
\begin{equation}\label{ratineq}
N \leq (148d\sigma^2)\log N + 971d\sigma^2.
\end{equation}
Assuming for now that this holds, it follows by Lemma~\ref{tlem} that
\begin{equation*}
\begin{split}
N & \leq \Big(\frac{e}{e-1}\Big)(148d\sigma^2\log(148d\sigma^2)+971d\sigma^2) \\
	 & = \Big(\frac{e}{e-1}\Big)148d\sigma^2\log(148e^{971/148}d\sigma^2).
\end{split}
\end{equation*}
In view of this and the fact that $|S|\leq 24^2N$, the bound $(\ref{kratmain})$ follows immediately. 

\medskip

It now remains only to prove $(\ref{ratineq})$.  Let $v_0$ be an archimedean place of $k$, chosen so that $|j_E|_{v_0}=\max_{v\mid\infty}(|j_E|_{v})$.  We then have a corresponding embedding $\sigma:E(k)\hookrightarrow\C/L$, where $L=\Z+\tau\Z$ is a normalized lattice, and $|j(\tau)|=|j_E|_{v_0}$.  If $\lambda:(\C/L)\setminus\{0\}\to\R$ denotes the N\'eron function on the complex torus, then $\lambda_{v_0}=\lambda\circ\sigma$. 

Divide the torus $\C/L$ into the $24^2$ parallelograms 
\begin{equation*}
P_{m_1,m_2} = \Big\{z=r_1+r_2\tau\,\Big|\, \frac{m_1-1}{24}\leq r_1\leq \frac{m_1}{24}\text{ and } \frac{m_2-1}{24}\leq r_2\leq \frac{m_2}{24}\Big\},
\end{equation*}
where $1\leq m_1,m_2\leq 24$.  By the pigeonhole principle there exists a set $Z=\{P_1,\dots,P_N\}\subseteq S$ of $N$ distinct points 
such that $\sigma(Z)$ is contained in one of the $24^2$ parallelograms.  In particular, 
it follows that the difference $\sigma(P_i)-\sigma(P_j)$ between any two points in $\sigma(Z)$ must lie in one of the four parallelograms 
$P_{m_1,m_2}$, $m_1,m_2\in\{1,24\}$.  Therefore, by Lemma~\ref{hind} we have $\lambda_{v_0}(P_i-P_j)\geq\frac{1}{288}\max\{1,\log|j_E|_{v_0}\}$ 
for all such pairs.  

It follows from the above considerations that
\begin{equation}\label{v0estimate}
\Lambda_{v_0}(Z) \geq \frac{(N-1)}{288N}\max\{1,\log|j_E|_{v_0}\};
\end{equation}
and by Lemma~\ref{lambdalemmaarch} and the maximality of $|j_E|_{v_0}$, we have
\begin{equation}\label{nonv0estimate}
\begin{split}
\sum_{\stackrel{v\mid\infty}{v\neq v_0}}\frac{d_v}{d}\Lambda_v(Z) & \geq -\sum_{\stackrel{v\mid\infty}{v\neq v_0}}\frac{d_v}{d} 
      \Big\{\frac{\log N}{2N} +\frac{1}{12N}\log^+|j_E|_{v} + \frac{16}{5N}\Big\} \\
        & \geq -\frac{\log N}{2N} - \frac{1}{12N}\log^+|j_E|_{v_0} - \frac{16}{5N}.
\end{split}
\end{equation}
If $N\leq24d$, then $(\ref{ratineq})$ plainly holds; so we may henceforth assume that $N\geq24d+1$.  Then combining the estimates $(\ref{v0estimate})$ and $(\ref{nonv0estimate})$ we have
\begin{equation}\label{ilb}
\begin{split}
\sum_{v\in\M_k^\infty}\frac{d_v}{d}\Lambda_v(Z) & \geq -\frac{\log N}{2N}-\frac{16}{5N} +\frac{1}{N}\Big\{\frac{(N-1)d_{v_0}}{288d}-\frac{1}{12}\Big\}\max\{1,\log|j_E|_{v_0}\} \\
        & \geq  -\frac{\log N}{2N}-\frac{16}{5N} +\frac{1}{N}\Big\{\frac{N-1}{288d}-\frac{1}{12}\Big\} \\
        & =  -\frac{\log N}{2N}-\frac{197}{60N} + \frac{(N-1)}{288dN}.
\end{split}
\end{equation}

\medskip

We now turn to the lower bounds on $\Lambda_v(Z)$ at the non-archimedean places; in particular we will show that
\begin{equation}\label{brlb}
\sum_{v\in\M_k^0}\frac{d_v}{d}\Lambda_v(Z) \geq  \frac{1}{12d}\Big(\frac{1}{16\sigma^2} - \frac{1}{N}\Big)\log |\N(\Delta)|.
\end{equation}
First, if $E/k$ has everywhere good reduction then $(\ref{brlb})$ plainly holds, since in that case $\log |\N(\Delta)|=0$ and the left-hand-side is nonnegative by Lemma~\ref{lambdalemma}.  So we may assume that the set $\M_k^{\br}=\{v\in\M_k^{0}\,\mid\, E/k_v \text{ has bad reduction} \}$ is nonempty.

At this point we will require Jensen's inequality: if $\{w_v\}$ is a finite set of positive weights with $\sum_v w_v=W$, if $x_v>0$ for all $v$, and if $\phi(x)$ is a convex function for $x>0$, then 
\begin{equation}
\sum_v w_v\phi(x_v) \geq W\phi\Big(\frac{1}{W}\sum_v w_v x_v\Big).
\end{equation}
Applying this with $\phi(x)=1/x$, with $\{w_v=\eta_v\log|\N(\p_v)| \,\mid\, v\in\M_k^{\br}\}$ as our weights (thus $W=\log|\N(\f)|$), and with $x_v=\eta_vc_v^2/\delta_v$, we have
\begin{equation}\label{ji}
\begin{split}
\sum_{v\in\M_k^{\br}} \frac{\delta_v}{c_v^2}\log|\N(\p_v)| & \geq W\Big(\frac{1}{W}\sum_{v\in\M_k^{\br}} \frac{\eta_v^2c_v^2}{\delta_v}\log|\N(\p_v)| \Big)^{-1} \\
	& \geq W^2\Big(\sum_{v\in\M_k^{\br}}16\delta_v\log|\N(\p_v)| \Big)^{-1} \\
	& = \frac{(\log\N(\f))^2}{16\log|\N(\Delta)|}.
\end{split}
\end{equation}
The first inequality in $(\ref{ji})$ is Jensen's, and the second follows from the inequality 
\begin{equation}\label{localbound}
\frac{\eta_v^2c_v^2}{\delta_v} \leq 16\delta_v
\end{equation}
for $v\in\M_k^{\br}$.  (To see this, recall that if $E/k_v$ has split multiplicative reduction then $\delta_v=c_v$ and $\eta_v=1$, and $(\ref{localbound})$ holds; while otherwise $\eta_v\leq\delta_v$ and $c_v\leq4$, and $(\ref{localbound})$ follows in this case as well.)

Finally, combining $(\ref{lambdalowerbound})$ and $(\ref{ji})$ and noting that $d_v\log|1/\Delta_v|_v=\delta_v\log|\N(\p_v)|$, we have
\begin{equation*}
\begin{split}
\sum_{v\in\M_k^{0}}\frac{d_v}{d}\Lambda_v(Z) & \geq  \sum_{v\in\M_k^{\br}}\frac{d_v}{d}\Big(\frac{1}{c_v^2}-\frac{1}{N}\Big)\frac{1}{12}\log|1/\Delta_v|_{v}  \\
	& = \frac{1}{12d}\sum_{v\in\M_k^{\br}}\Big(\frac{1}{c_v^2}-\frac{1}{N}\Big)\delta_v\log|\N(\p_v)|  \\
	& \geq \frac{1}{12d}\Big(\frac{(\log|\N(\f)|)^2}{16\log|\N(\Delta)|}-\frac{1}{N}\log|\N(\Delta)|\Big)  \\
	& = \frac{1}{12d}\Big(\frac{1}{16\sigma^2} - \frac{1}{N}\Big)\log |\N(\Delta)|,
\end{split}
\end{equation*}
which is $(\ref{brlb})$.

\medskip

We are now ready to combine these local estimates.  By the parallelogram law we have the global upper bound 
\begin{equation}
\begin{split}
\Lambda(Z) & \leq \frac{4}{N}\sum_{j=1}^{N}\h(P_j) \\
        & \leq \frac{\log |\N(\Delta)|}{2^{11}3d\sigma^2}
\end{split}
\end{equation}
on the sum defined in $(\ref{heightdisc})$, by the upper bound on $\h(P)$ for $P\in S$.  Therefore by $(\ref{ilb})$ and $(\ref{brlb})$ we have
\begin{equation}\label{genglobal}
\begin{split}
\frac{\log |\N(\Delta)|}{2^{11}3d\sigma^2} & \geq \Lambda(Z) \\
	& = \sum_{v\in\M_k}\frac{d_v}{d}\Lambda_v(Z) \\
        & \geq  -\frac{\log N}{2N}-\frac{197}{60N} + \frac{(N-1)}{288dN} + \frac{1}{12d}\Big(\frac{1}{16\sigma^2} - \frac{1}{N}\Big)\log |\N(\Delta)|.
\end{split}
\end{equation}

If $E/k$ has everywhere good reduction then, then $\sigma=1$ and $\log |\N(\Delta)|=0$, and thus $(\ref{genglobal})$ becomes 
\begin{equation}\label{goodglobal}
0  \geq  -\frac{\log N}{2N}-\frac{197}{60N} + \frac{(N-1)}{288dN}.
\end{equation}
It follows that $N\leq 144d\log N + 945.6d+1 \leq 144d\log N + 946.6d$, and $(\ref{ratineq})$ holds in this case.  On the other hand, suppose that $\M_k^{\br}$ is nonempty.  If $N\leq2^9\sigma^2$, then $(\ref{ratineq})$ holds, so we may assume that $N>2^9\sigma^2$.  The bound $(\ref{genglobal})$ implies that 
\begin{equation}\label{badglobal}
\begin{split}
0 & \geq  -\frac{\log N}{2N}-\frac{197}{60N} + \frac{1}{12d}\Big(\frac{1}{16\sigma^2} - \frac{1}{N}\Big)\log |\N(\Delta)| -\frac{\log |\N(\Delta)|}{2^{11}3d\sigma^2} \\
        & \geq -\frac{\log N}{2N}-\frac{197}{60N} + \frac{1}{12d}\Big(\frac{1}{2^4\sigma^2} - \frac{1}{2^9\sigma^2}-\frac{1}{2^9\sigma^2}\Big)\log |\N(\Delta)| \\
        & = -\frac{\log N}{2N}-\frac{197}{60N} + \frac{5}{2^{10}d\sigma^2}\log |\N(\Delta)| \\
        & \geq -\frac{\log N}{2N}-\frac{197}{60N} + \frac{5\log 2}{2^{10}d\sigma^2},
\end{split}
\end{equation}
since $\log |\N(\Delta)|\geq\log2$.  We deduce that 
\begin{equation}
\begin{split}
N & \leq \Big\{\frac{2^9d\sigma^2}{5\log2}\Big\}\log N + \frac{197\cdot2^{10}d\sigma^2}{60\cdot5\log2} \\
	& < (148d\sigma^2)\log N + 971d\sigma^2.
\end{split}
\end{equation}
Thus $(\ref{ratineq})$ holds in this case as well, and the proof of Proposition~\ref{smallprop} is complete.
\end{proof}

\begin{proof}[Proof of Theorems \ref{krationaltor} and \ref{langs}]
Again, let $S$ denote the set on the left-hand-side of $(\ref{kratmain})$.  Theorem~\ref{krationaltor} follows trivially from Proposition~\ref{smallprop}, since $E(k)_{\tor}\subseteq S$.  To see that Theorem~\ref{langs} follows, let $P\in E(k)$ be a non-torsion point, and let $M$ be the largest integer such that $\h((M-1)P)\leq (\log |\N_{k/\Q}(\Delta_{E/k})|)/2^{13}3d\sigma^2$.  Then the first $M$ multiples $O, P, 2P, \dots, (M-1)P$ of $P$ are contained in the set $S$, and therefore $M\leq c_1d\sigma^2\log(c_2d\sigma^2)$ by Proposition~\ref{smallprop}.  On the other hand, by the maximality of $M$ we have $M^2\h(P)=\h(MP)>(\log |\N_{k/\Q}(\Delta_{E/k})|)/2^{13}3d\sigma^2$.  
Therefore 
\begin{equation*}
\begin{split}
\h(P) & > M^{-2}\frac{\log|\N_{k/\Q}(\Delta_{E/k})|}{2^{13}3d\sigma^2} \\
        & \geq c(d,\sigma)\log|\N_{k/\Q}(\Delta_{E/k})|,
\end{split}
\end{equation*}
where $c(d,\sigma)$ is given by $(\ref{cdsigma})$.
\end{proof}
\end{section}


\end{document}